\documentclass[12pt]{amsart}

\usepackage{amssymb}
\usepackage{latexsym}
\usepackage{verbatim}
\usepackage{fullpage}

\input {cyracc.def}
\font\ququ=cmr10 scaled \magstep1
\font\tencyr=wncyr8 scaled \magstep1
\def\rus{\tencyr\cyracc}

\newcommand{\re}[1]{\textrm  (\ref{#1})}

\renewenvironment{proof}
{\noindent {\sl Proof.}\quad }{\hfill q.e.d.
\vskip1.1ex\noindent }

\newenvironment{proof*}
{\noindent {\sl Proof.}\quad }{\hfill 
q.e.d.}

\renewcommand{\theequation}{\thesection$\cdot$\arabic{equation}}
\renewcommand{\thesubsubsection}{\theequation$\cdot$\arabic{subsubsection}}

\catcode`\@=\active
\catcode`\@=11
\def\@eqnnum{\hbox to
.01pt{}\rlap{\hskip-\displaywidth(\mathbf{\theequation})}}

\renewcommand{\@cite}[2]{[{{\bf #1}\if@tempswa , #2\fi}]}
\renewcommand{\@biblabel}[1]{[{\bf #1}]\hfill}

\catcode`\@=12
\newenvironment{s}[1]
{ \vskip1.2ex \refstepcounter{equation}
\noindent {\bf \theequation\quad #1.} \begin{sl}}{\end{sl}
\vskip1.1ex\noindent }

\newenvironment{rem}[1]
{ \vskip1.2ex \refstepcounter{equation}
\noindent {\bf \theequation\enspace {#1}.} }{ \vskip1.1ex\noindent }

\newenvironment{rems}[1]
{ \vskip1.2ex \refstepcounter{subsubsection}
\noindent {\bf \thesubsubsection\enspace {#1}.} }{ \vskip1.1ex\noindent }

\newenvironment{subs}[1]
{\vskip1.2ex \refstepcounter{equation}
\noindent {\bf (\theequation)\quad #1.} }{\quad}

\newcommand {\sekt}[1]
{{\vskip2.5ex\refstepcounter{section}\setcounter{equation}{0}
\setcounter{subsubsection}{0}
\noindent\large \bf \thesection\quad \parbox[t]{424pt}{#1}
\nopagebreak\vskip1.5ex\noindent}}

\newcommand {\ah}{{\frak a}}
\newcommand {\be}{{\frak b}}
\newcommand {\de}{{\frak d}}

\newcommand {\g}{{\frak g}}
\newcommand {\h}{{\frak h}}
\newcommand {\ka}{{\frak k}}
\newcommand {\el}{{\frak l}}
\newcommand {\n}{{\frak n}}

\newcommand {\p}{{\frak p}}
\newcommand {\q}{{\frak q}}
\newcommand {\rr}{{\frak r}}
\newcommand {\es}{{\frak s}}
\newcommand {\te}{{\frak t}}
\newcommand {\ut}{{\frak u}}
\newcommand {\z}{{\frak z}}

\newcommand {\glv}{{\frak gl}(V)}

\newcommand {\ap}{\alpha}

\newcommand {\lb}{\lambda}

\newcommand {\vp}{\varphi}
\newcommand {\ca}{{\mathcal A}}
\newcommand {\cb}{{\mathcal B}}
\newcommand {\ck}{{\mathcal K}}
\newcommand {\N}{{\mathcal N}}

\newcommand {\md}{/\!\!/}

\newcommand {\isom}{\stackrel{\sim}{\longrightarrow}}
\newcommand {\ad}{{\mathrm{ad\,}}}

\newcommand {\hot}{{\mathrm{ht\,}}}
\newcommand {\ind}{{\mathrm{ind\,}}}

\newcommand {\Ker}{{\mathrm{Ker\,}}}
\newcommand {\Ima}{{\mathrm{Im\,}}}
\newcommand {\rank}{{\mathrm{rank\,}}}
\newcommand {\rk}{{\mathrm{rk\,}}}
\newcommand {\spe}{{\mathrm{Spec\,}}}

\newcommand {\tri}{{\frak sl}_2}

\newcommand {\ov}{\overline}
\newcommand {\un}{\underline}

\newcommand {\vno}[1]{\vskip#1 ex\noindent}
\newcommand {\rar}{\rightarrow}
\newcommand {\qus}{\hfill $\square$ \vno{1.1}}
\newcommand {\qu}{\hfill $\square$}
\newcommand {\beq}{\begin{equation}}
\newcommand {\eeq}{\end{equation}}

\newcommand{\nge}{\n_\g(e)}
\newcommand{\zge}{\z_\g(e)}
\newcommand{\zgc}{\z_\g(c)}
\newcommand{\zghf}{\z_\g(\hat f)}
\newcommand{\dge}{\de_\g(e)}

\newfam\Bbbfam\newfam\eufam%
\font\Bbbfont=msbm10 scaled 1200%
\font\olala=msam10 scaled 1200%
\font\frak=eufm10 scaled 1400%
\font\Bbbsmallfont=msbm8%
\textfont\Bbbfam=\Bbbfont\scriptfont\Bbbfam=\Bbbsmallfont%%
\font\euzw=eufm10 scaled 1200%
\font\euac=eufm7 scaled 1200%
\font\euacc=eufm7 scaled 1000%
\textfont\eufam=\euzw\scriptfont\eufam=\euac%
\scriptscriptfont\eufam=\euacc%
\def\frak{\fam\eufam}%
\def\Bbb{\fam\Bbbfam}%
\def\bbb{\fam\Bbbfam\scriptstyle}
\def\varnothing{\hbox {\Bbbfont\char'077}}
\def\square{\hbox {\olala\char"03}}
\def\bbk{\hbox {\Bbbfont\char'174}}
\def\semidir{\hbox {\olala\char"68}}

\begin{document}
\setlength{\parskip}{2pt plus 4pt minus 0pt}
%\hfill {\scriptsize \today} \vskip1ex
\vskip1ex

\begin{center}
{\Large \bf %\parbox[t]{\textwidth}
{The index of a Lie algebra, 
the centraliser of a nilpotent element, and the normaliser of the
centraliser}} 
\medskip \\
{By DMITRI I. PANYUSHEV}\medskip  \\%
%\footnote{This research was supported in part by 
%CRDF Grant No.\, RM1--2088}
%\hfill {\scriptsize preliminary version}
{\footnotesize 
%\parbox{195pt}{% 
{\it Independent University of Moscow,
Bol'shoi Vlasevskii per. 11 \\ 
121002 Moscow, \quad Russia \\ e-mail}: {\tt panyush@mccme.ru }
%\\{(\it Received 30 April 2001})}
}
\end{center}

\medskip
\smallskip

\noindent {\large \bf Introduction}
\vno{2}%
Let $\q$ be a Lie algebra and $\xi$ a linear form on $\q$. Denote by
$\q_\xi$ the set of all $x\in \q$ such that $\xi([\q,x])=0$. In other
words, $\q_\xi=\{x\in\q\mid (\ad^*x){\cdot}\xi=0\}$, where 
$\ad^*:\q\to {\frak gl}(\q^*)$ is the coadjoint representation of $\q$.
%$\q^\ast$ its dual space. Then
%$\q$ acts on $\q^\ast$ via the coadjoint representation.
The {\it index\/} of $\q$, $\ind \q$, is defined by
\[ 
   \ind\q =\min_{\xi\in\q^*}\dim\q_\xi \ .
\]
This definition goes back to J.\,Dixmier, see \cite[11.1.6]{Di}. He considered
index because of its importance in Representation Theory.
The problem of computing the index may also
be treated as part of Invariant Theory.
For, if $\q$ is an algebraic Lie algebra and $Q$ is a corresponding 
algebraic group,
then $\ind\q$ equals the transcendence degree of the field of $Q$-invariant
rational functions on $\q^*$.
If $\q$ is reductive, then $\q$ and $\q^\ast$ are isomorphic as 
$\q$-modules and hence $\ind\q=\rk\q$. It is therefore interesting to
study index for non-reductive Lie algebras. 
On the other hand, studying the index for all Lie algebras is too
pretentious. Therefore I think that the most promising approach is to
look at the index for various natural classes of non-reductive subalgebras of
semisimple Lie algebras. There are at least two such classes:
\ a) parabolic subalgebras and their `relatives` (nilpotent radicals, 
seaweeds, etc.); \ b) centralisers of elements and their `relatives'.
Some recent results on a) are found in \cite{mmj}, whereas the present paper
deals mainly with b).
\\[.6ex]
The ground field $\bbk$ is algebraically closed and of characteristic 0.
Let $\g$ be a semisimple Lie algebra and $G$ its adjoint group.
Given $x\in\g$, let $\z_\g(x)$ denote the centraliser of $x$. The normaliser
in $\g$ of $\z_\g(x)$ is denoted by $\n_\g(x)$.
The goal of this paper is to study $\ind\z_\g(x)$ and $\ind\n_\g(x)$
for $x$ nilpotent.  By a result of R.\,Brylinski and B.\,Kostant
\cite{rk1}, $\n_\g(x)$ is determined by $\z_\g(x)$ and 
%%the centre of $\z_\g(x)$ or 
the double centraliser of $x$, denoted by $\de_\g(x)$. 
For this reason, we include 
a discussion of some properties of $\de_\g(x)$.  We also prove some general 
results on the index, which are of independent interest. The following brief 
exposition of our results is not exhaustive. 
\\[.6ex]
In Section 1, we collect some old and new general results for the index
of an arbitrary Lie algebra. Actually, to deal with index of Lie algebras,
a more general concept of the index of a representation is required. 
Using this,
we prove an inequality relating the index of an algebra and an ideal in it.
More precisely, let $\q$ be an ideal in $\tilde\q$. Then 
\[
   \ind\q+\ind\tilde\q\le \dim(\tilde\q/\q)+ 2\,\ind(\tilde\q,\q) \ ,
\]
where $\ind(\tilde\q,\q)$ is the index of $\q$ as $\tilde\q$-module. In 
particular, $\p^u$ being the nilpotent radical of a parabolic subalgebra
$\p\subset\g$, we obtain $\ind\p + \ind\p^u\le \dim (\p/\p^u)$.
\\[.6ex]
Section 2 begins with generalities on $\tri$-triples and
centralisers of nilpotent elements. Some properties of $\de_\g(x)$ are
discussed afterwards. Let $\{x,h,y\}$ be an $\tri$-triple. 
As is known, $h$ determines a $\Bbb Z$-grading of $\g$, which is
nonnegative 
on $\z_\g(x)$ and hence on $\de_\g(x)$. 
We give a simple proof for the result of
R.\,Brylinski and Kostant \cite{rk2}
that $\de_\g(x)(2)=\bbk x$, if $\g$ is simple.
Another result is that the Killing form is non-degenerate on
$\de_\g(x)\oplus\de_\g(y)$. This is equivalent to the fact that
$\de_\g(y)\oplus [\g,\z_\g(x)]=\g$.
\\[.6ex]
In Section 3, we discuss a conjecture of Elashvili to the effect that
$\ind\z_\g(x)=\rk\g$ for all $x\in\g$. The Jordan decomposition immediately
reduces this
conjecture to nilpotent elements. We show that the conjecture holds for
`large' and `small' nilpotent $G$-orbits. The claim for the regular
and subregular orbits follows 
from the fact that $\z_\g(x)$ is Abelian if and only if $x$ regular.
The case-by-case proof of the latter was known, 
but we give a short conceptual proof, see Theorem~\ref{concept}.
This solves Problem~15 posed by Steinberg at the Moscow I.C.M., 
see~\cite{icm66}. Our proof exploits work of Springer \cite{tony} and
the above-mentioned result of Brylinski and Kostant.
At the other extreme, we prove that Elashvili's conjecture
holds for $x$ such that $(\ad x)^3=0$.
\\[.6ex]
In Section 4, we consider $\n_\g(x)$ and its natural representation on
$\de_\g(x)$. It is proved that $\ind\n_\g(x)=\ind\z_\g(x)-\dim\de_\g(x)$
whenever there exists a $\xi\in\de_\g(x)^*$ such that
$\n_\g(x){\cdot}\xi=\de_\g(x)^*$.
%%$N_\g(x)$ has a dense orbit in $\de_\g(x)^*$. Here $N_\g(x)$
%%is the connected algebraic group with Lie algebra $\n_\g(x)$.
We also show that the last condition (and even a stronger one)
is satisfied for all $x$ if
$\g\in\{ {\frak sl}_{n}, {\frak sp}_{2n}, {\frak so}_{2n+1}\}$
and for $x$ in $\g={\frak so}_{2n}$ such that the corresponding partition
has at least three parts, see
Theorem~\ref{3serii}.
\\[.6ex]
In Section 5, we prove that if $x$ is regular, then $\n_\g(x)$ is a 
{\it Frobenius\/}
Lie algebra, i.e., it is of index 0. It is worth mentioning that each
Frobenius subalgebra of $\g$ determines a solution of the
classical Yang-Baxter equation (CYBE), see \cite[\S\,8]{bd}. 
Hence our result means that the
regular nilpotent orbit yields a solution of CYBE.
\\[.6ex]
We conjecture that there is a simple general formula for the index of
$\n_\g(x)$, namely
$\ind\n_\g(x)=\rk\g-\dim\de_\g(x)$. Our results show that it is true
if either $x$ is regular or $(\ad x)^3=0$. 
%Furthermore, modulo 
%Elashvili's conjecture, the formula for $\ind\n_\g(x)$ holds for
%all $x$ in 
This can be thought of as more general 
version of Elashvili's conjecture, as the conjecture for $\ind\n_\g(x)$
implies that for $\ind\z_\g(x)$.
\\[.6ex]
{\footnotesize {\it Acknowledgements.} I am grateful to A.G.\,Elashvili 
for sharing with me his enthusiasm about the index of centralisers. 
%%and for numerous conversations on this subject. 
Thanks are also due to
P.\,Tauvel for interesting discussions on the double centralisers.
This research was supported in part by Alexander von Humboldt-Stiftung, 
RFBI Grant 01--01--00756, and
CRDF Grant {RM1-2088}. }

\sekt{The index of a Lie algebra
\label{index}}%
The index of a Lie algebra $\q$, as defined in the Introduction, is
%%may be regarded as 
an integer attached to the adjoint representation of $\q$. 
A similar procedure applied to an arbitrary representation
of $\q$ yields the notion of the {\it index of a representation\/}.
This is quite useful, even though one is mainly interested in the
index of Lie algebras. Let $\rho:\q \rar \glv$ be a 
finite-dimensional representation of $\q$, i.e. $V$ is a $\q$-module.
Abusing notation, we write $s{\cdot}v$ in place of $\rho(s)v$, if $s\in\q$
and $v\in V$. An element $v\in V$ 
is called {\it regular\/} or $\q$-{\it regular\/} whenever its
stabiliser $\q_v=\{s\in\q\mid s{\cdot}v=0\}$ has minimal dimension. Because
the function  $v \mapsto \dim{\q}_v$ ($v\in V$) is upper semicontinuous,
the set of all $\q$-regular elements is open and dense in $V$. 
Consider also the dual $\q$-module $V^\ast$. 
\\[1ex]
{\bf Definition.} 
The nonnegative integer
$\displaystyle\dim V - \max_{\xi\in V^*}(\dim \q{\cdot}\xi)$ 
is called the {\it index\/} of (the $\q$-module) $V$. 
It will be denoted by $\ind(\rho, V)$ or $\ind (\q, V)$.
\\[.8ex]
Notice that 
in order to define the index of $V$ we used the regular elements in $V^*$\,!
It is clear that $\ind(\ad,\q)=\ind\q$ is the index of $\q$ 
in the sense of Introduction.  
%%For future reference, we collect some  easy properties of this index.
\\[.6ex]
Consider the bilinear form with values in $V$
\[
\ck=\ck(\q,V) : \q\times V \to V \ ; \quad (s,v)\mapsto s{\cdot}v
\]
Evaluating this form against an arbitrary
element $\xi\in V^\ast$ gives a form with values in $\bbk$:
\[
\ck_\xi: \q\times V\to V\stackrel{\xi(\cdot)}{\longrightarrow} \Bbbk \ .
\]
We may consider $\ck_\xi$ as an element of
$\mathrm{Hom} (\q,V^*)$. The following is obvious.
\begin{s}{Lemma} \label{refer} 
$\Ker (\ck_\xi)=\q_\xi$ and $\Ima (\ck_\xi)=\q{\cdot}\xi$.
\end{s}%
It follows that $\ind(\q,V)=\dim V - \displaystyle
\max_{\xi\in V^*} (\rank \ck_\xi)$.
Let $n=\dim\q$ and $m=\dim V$. Having chosen bases for $\q$ and $V$,
we may regard $\ck$ as $n\times m$-matrix with entries in $V$, where
$V$ is identified with the component of grade 1 in the symmetric algebra 
$\mathcal S  ^\bullet(V)$. Then
$\rank\ck=\max_{\xi\in V^*}(\rank\ck_\xi)$. Thus,
\begin{equation} \label{ind(q,V)}
\ind(\q,V)=\dim V- \rank\ck(\q,V) \ .
\end{equation}
In case $V=\q$, we see that $\ck=\ck(\q)$ is nothing but the Lie bracket, 
i.e., $\ck(s',s'')=[s',s'']$, $s',s''\in\q$.
Notice also that
%\begin{s}{Lemma} \label{reference} \\
%1. If $\q=\q_1\dotplus\q_2$, %% is a direct sum of Lie algebras, 
%then $\ind\q=\ind\q_1+\ind\q_2$; \\ 
%%$\ind\q=\dim\q$ if and only if $\q$ is Abelian; \\
%2. $\ind\q \ge \dim \z(\q)$, where $\z(\q)$ is the centre of $\q$;\\
%3. If $\q$ is reductive, then $\ind\q=\rk\q$. \qu
%\end{s}%
%
%\quad $K=K(\q)
%: \q\times\q\to\q$. Evaluating this form against an arbitrary
%point $\xi\in\q^\ast$ gives a usual skew-symmetric bilinear form:
\[
\ck(\q)_\xi: \q\times\q\to\q\stackrel{\xi(\cdot)}{\longrightarrow} \Bbbk \ .
\]
is the {\it Kirillov form\/} associated with $\xi\in\q^*$.
Since $\rank \ck(\q)_\xi$ is even, we deduce that
$\dim\q-\ind\q$ is even for any Lie algebra. 
\\
A Lie algebra $\q$
is called {\it Frobenius\/}, if $\ind\q=0$, i.e., $\ck(\q)_\xi$ is 
nonsingular for some $\xi\in\q^*$. 
\\[.7ex]
In case $\q$ is an algebraic Lie algebra, a more geometric description
is available. Denote by $Q$ an algebraic group with Lie algebra $\q$.
Then $\ind\q=\dim\q-\max_{\xi\in\q^*}\dim Q{\cdot}\xi$. In particular,
$\q$ is Frobenius if and only if $Q$ has an open orbit in $\q^*$.
\\[.7ex]
A useful tool for computing index is a theorem of M.\,Ra\"\i s.
It also shows that the index of a representation is helpful in computing the
index of Lie algebras.
Let $\rho:\q \rar \glv$ be any representation of $\q$.
The linear space $V\times \q$ has a natural structure of Lie algebra,
with bracket $[\ ,\ ]\widetilde{\ }$ defined by the equality
\[
[(v_1,s_1),(v_2,s_2)]\widetilde{\ }=(\rho(s_1)v_2-\rho(s_2)v_1, [s_1,s_2])\ .
\]
Following Ra\"\i s, the resulting Lie algebra is denoted by 
$V\times_\rho \q$. It is 
a semi-direct product of $\q$ and an Abelian ideal $V$. We identify the
dual space $(V\times_\rho\q)^*$ with $V^*\oplus\q^*$.
\begin{s}{Theorem {\ququ (Ra\"\i s \cite{rais})}}  \label{thm:rais}
Let $\xi\in V^\ast$ be a $\q$-regular element such that $\xi+\eta$ is a 
regular element in $(V\times_\rho\q)^*$ for some $\eta\in\q^*$. Then
$\ind(V\times_\rho\q)= \ind(\q, V)+\ind\q_\xi$.
\end{s}%
It is not assumed %in the Theorem 
that the stabilisers of all (or almost all) regular elements are 
conjugated. 
%%We are going to apply Theorem~\ref{thm:rais} in the following 
%%specific situation.
%\begin{s}{Corollary}  \label{cor:rais}
%Suppose $\q$ is an algebraic Lie algebra and
%$\ind(\q:V)=0$. Then 
%$\ind(V\times_\rho\q)= \ind\q_\xi$, where $\xi\in V^*$ is \un{any} 
%$\q$-regular element.
%\end{s}\begin{proof}
%Let $Q$ be a connected algebraic group with Lie algebra $\q$. The condition
%on  $\ind(\q:V)$ is equivalent to the fact that $Q$ 
%has an open orbit in $V^*$. Therefore the set of $\q$-regular
%elements in $V^*$ is just this orbit.
%\end{proof}%
%
%
The following result, which relates the indices
of a Lie algebra and an ideal in it, 
appears to be helpful in various applications.

\begin{s}{Theorem}   \label{atilde}  \\
Let $\q$ be an ideal in a Lie algebra $\tilde\q$.
Then $\ind\q+\ind\tilde\q\le \dim(\tilde\q/\q)+ 2\,\ind(\tilde\q,\q)$.  
\end{s}\begin{proof}
Let $m=\dim\q$ and $n=\dim\tilde\q$, $m \le n$. 
Choose a basis $\{e_1,\dots,e_n\}$ for $\tilde\q$ so that
$\{e_1,\dots,e_m\}$ is a basis for $\q$. 
The $n\times n$-matrix $\ck(\tilde\q)=([e_i,e_j])_{i,j=1}^n$ 
has the following block structure
in this basis:
\[
  \ck(\tilde\q)=\left(\begin{array}{cr} C & B \\ -B^t &  D
\end{array}\right) \ .
\]
Here $C=\ck(\q)$ is a skew-symmetric
matrix of order $m$ and $B$ is a rectangular 
$m\times(n-m)$-matrix.
Consider the block matrix
of order $m+n$ that depends on a parameter $\mu\in\bbk$:  
\[
M(\mu)=\left(\begin{array}{ccc} C    &   B    & C        \\ 
                                -B^t & \mu D  & -\mu B^t \\
                                C    & \mu B  &  \mu C
\end{array}\right) \ .
\]
We look at the rank of $M(0)$ in two different ways, which gives us
the required inequality. 

1) Since $(B\ C)^t=\left(\begin{array}{c} B^t \\ -C 
                   \end{array}\right)$, we have
$\textrm{rank\,}M(0)\ge 2\,\textrm{rank\,}(B\ C)$.
Also, up to a permutation of columns, the rectangular matrix
$(B\ C)$ is nothing but $\ck(\tilde\q,\q)$. Thus,
$\textrm{rank\,}M(0)\ge \displaystyle 2\max_{\xi\in\q^*}\dim
(\tilde\q{\cdot}\xi)=2(\dim\q-\ind(\tilde\q,\q))$. 

2) If $\mu\ne 0$, then a chain of elementary transformations brings $M(\mu)$
in the following form: 
\[
  M(\mu) \sim 
       \left(\begin{array}{ccc} C      &   B    & (1-\mu)C        \\ 
                                -B^t   & \mu D  &  0 \\
                             (1-\mu)C  &   0    &  (\mu-1)\mu C
             \end{array}\right)
         \sim
 \left(\begin{array}{ccc} \frac{1}{\mu}C  &   B    &  0        \\ 
                                   -B^t   & \mu D  &  0 \\
                                      0   &   0    &  (\mu-1)\mu C
             \end{array}\right) ß .
\]
Therefore if $\mu\ne 0,1$, then 
\[\rank M(\mu)=\rank C +
\rank \left(\begin{array}{cr} \frac{1}{\mu}C & B \\ -B^t &  \mu D
                 \end{array}\right) = 
                 \rank C +
\rank \left(\begin{array}{cr} C & B \\ -B^t &  D
                 \end{array}\right) \ .
                 \]
Thus,
\[
\rank M(\mu)=\rank \ck(\q) + \rank \ck(\tilde\q)=\dim\q-\ind\q+\dim\tilde\q
-\ind\tilde\q \ .
\] 
From the semi-continuity of the rank, it follows that 
$\rank M(\mu)\ge \rank M(0)$ for $\mu\ne 0,1$, which yields the 
required inequality.                
\end{proof}%
Let us apply this theorem to parabolic subalgebras. Let $\p$ be a
%%the case, where  $\tilde\q=\p$ is a 
parabolic subalgebra
of a semisimple Lie algebra $\g$, with the nilpotent radical 
$\p^u$ and a Levi subalgebra $\el$. Let $\be\subset \p$ be any Borel 
subalgebra of $\g$. Obviously, $\be\supset\p^u$ and $\p^u$ is a $\be$-module.
It is known that the Borel subgroup $B\subset G$ has a dense
orbit in $(\p^u)^*$, i.e., $\ind(\be,\p^u)=0$, see \cite[2.4]{joseph}. This 
clearly implies that $\ind(\p,\p^u)=0$.
\begin{s}{Corollary}  \label{levi} 
For any parabolic subalgebra $\p$ of $\g$ we have 

{\sf (i)} \ $\ind\p +\ind\p^u \le \dim\el$.
%%where $\el$ is a Levi subalgebra of $\p$.
In particular, $\ind\be +\ind\be^u \le \rk\g$.

{\sf (ii)} \ $\ind\be +\ind\p^u \le \dim (\be/\p^u)$.
%%where 
\end{s}\begin{proof}
Take $\q=\p^u$ and then $\tilde\q=\p$ or $\be$.
\end{proof}%
\vskip-1.5ex
\begin{rems}{Remarks} Note that $\dim (\be/\p^u)$ is dimension of a Borel
subalgebra of $\el$.
It can conceptually be proved that $\ind\be^u=\rk\g^\theta$
and $\ind\be= \rk\g-\rk\g^\theta$, where $\theta$ is 
an involutory automorphism of $\g$ of maximal rank. In particular, 
$\ind\be +\ind\be^u =\rk\g$.
\end{rems}%
The next Proposition was communicated to me by A.G.\,Elashvili at the
beginning of 90's. He attributed the result to E.B.\,Vinberg.
Being rather simple, this result is extremely helpful.
As Vinberg's proof is not published, we give a proof below. 
\\
Let $\rho:\q \rar \glv$ be a representation of $\q$. Given $w\in V$,
consider the natural representation of $\q_w$ on $\ov{V}=V/\q{\cdot}w$.
\begin{s}{Proposition {\ququ (Vinberg)}}  \label{vinberg} 
For any $w\in V$, we have \\[.6ex]
\hbox to \textwidth{\hfil
$\displaystyle \max_{v\in V}\dim(\q{\cdot}v)\ge
         \max_{\eta\in \ov{V}}\dim(\q_w{\cdot}\eta) + \dim(\q{\cdot}w) \ .
$\hfil}\vskip-1ex
\end{s}\begin{proof}
Fix a vector space direct sum  $\q=\q_w\oplus\tilde\q$.
Choose $\eta=v+\q{\cdot}w\in\ov{V}$ so that $\dim(\q_w{\cdot}\eta)$
is maximal. For any $t\in \bbk\setminus\{0\}$, one has
$\q{\cdot}(tv+w)=\q_w{\cdot}v+ \tilde\q{\cdot}(tv+w)$. If $t=0$, then
dimension of the RHS is equal to $\dim(\q_w{\cdot}\eta)+\dim(\q{\cdot}w)$.
Consequently, for all but finitely many $t$'s, we obtain
$\dim\q{\cdot}(tv+w)\ge \dim(\q_w{\cdot}\eta)+\dim(\q{\cdot}w)$.
\end{proof}%
{\bf Remark.} If $\q_w$ is reductive, then one can identify $\ov{V}$ 
with a $\q_w$-stable complement to $\q{\cdot}w$. In this case, one can actually
prove that equality holds. 
It is an intersting open problem to find out weaker
sufficient conditions for the equality in \ref{vinberg} to hold.
See e.g. Conjecture~\ref{AG} below.
\begin{s}{Corollary}   \label{neravenstvo}
Let $\q$ be a Lie algebra. Then $\ind\q_\xi\ge
\ind\q$ for any $\xi\in\q^*$.
\end{s}\begin{proof*}
Apply the previous Proposition to $V=\q^*$. Here $\q^*/\q{\cdot}w\simeq
(\q_w)^*$.
\end{proof*}%
\begin{s}{Corollary {\ququ (Duflo--Vergne, 1969)}}   \label{duflo}
If $w\in\q^*$ is regular, then $\q_w$ is Abelian.
\end{s}\begin{proof}
Since $w$ is regular, $\ind\q=\dim\q_w$. Hence
$\dim\q_w\ge \ind\q_w\ge\ind\q=\dim\q_w$.
\end{proof}%

\sekt{Some properties of the double centraliser of a nilpotent element
\label{RK-prosto}}%
Let $\g$ be a semisimple Lie algebra with a fixed triangular decomposition 
$\g=\ut_-\oplus\te\oplus\ut_+$, $\Delta$ the corresponding root system, 
and $\Pi=\{\ap_1,\dots,\ap_p\}$ the set of simple roots. The Killing form on 
$\g$ is denoted by $\Phi$.
For any subset $M\subset\g$, let $Z_G(M)$ and $\z_\g(M)$ denote its
centraliser in $G$ and $\g$, respectively; $M^\perp$ is the orthogonal 
complement to $M$ with respect to the Killing form.
Let $\N \subset\g$ be
the nilpotent cone. By the Morozov-Jacobson theorem, each nonzero
element $e\in\N$ can be
included in an ${\frak sl}_2$-triple $\{e,h,f\}$ (i.e.,
$[e,f]=h,\ [h,e]=2e,\ [h,f]=-2f$). The semisimple element $h$
determines a $\Bbb Z$-grading in $\g$:
\[
 \g=\bigoplus_{i\in\Bbb Z}\g(i) \ ,
\]
where $\g(i)=\{\,x\in\g\mid [h,x]=ix\,\}$. Since
all ${\frak sl}_2$-triples containing $e$
are $Z_G(e)$-conjugate,
the properties of this $\Bbb Z$-grading 
does not depend on a particular choice of $h$.
\\[.5ex]
%Following E.B.\,Dynkin, we shall say that $h$ is a {\it characteristic\/}
%of $e$. 
The adjoint orbit $G{\cdot}h$ contains a unique element $h_+$ such that 
$h_+\in\frak t$ and $\alpha(h_+)\ge0$
for all $\alpha\in\Pi$.  
The Dynkin diagram of $\g$ equipped with the numerical
labels $\alpha_i(h_+)$, $\ap_i\in\Pi$, at the corresponding nodes 
is called the
{\it weighted Dynkin diagram} of $e$. After Dynkin, it is known 
(see  \cite[8.1,\,8.3]{dynkin}) that \par
(a) $\ap_i(h_+)\in\{0,1,2\}$; \par
(b) $\tri$-triples $\{e,h,f\}$ and $\{e',h',f'\}$
are $G$-conjugate {\sl if and only if\/}
%their characteristics 
$h$ and $h'$ are $G$-conjugate {\sl if and only if\/}
their weighted Dynkin diagrams coincide.
\vno{0.5}%
%\end{comment}
We also need the following standard results on the
structure of $Z_G(e)\subset G$ and 
$\z_\g(e)\subset\g$ (see e.g.~\cite[ch.\,III]{ss} or \cite{CoMc}). 
\begin{s}{Proposition}         \label{stab} 
\par
{\sf (i)} the Lie algebra $\z_\g(e)$ (resp.\ $\z_\g(f)$) is
positively (resp. negatively) graded; $\z_\g(e)=
\bigoplus_{i\ge 0} \z_\g(e)(i)$, where $\z_\g(e)(i)=\z_\g(e)
\cap\g(i)$, and likewise for $\z_\g(f)$; 
\par
{\sf (ii)} let $G(0)$ be the connected subgroup of $G$ with Lie
algebra $\g(0)$ and $K:= Z_{G(0)}(e)$ ({$K$ can be disconnected !}).
Then $K=Z_G(e)\cap Z_G(f)$ and it is a maximal reductive subgroup in 
both $Z_G(e)$
and $Z_G(f)$; \par
{\sf (iii)} for any $i$, there are $K$-stable
decompositions:
%*
\[ \g(i)=\z_\g(e)(i) \oplus[f,\g(i+2)],\qquad
   \g(i)=\z_\g(f)(i) \oplus[e,\g(i-2)].
\]
In particular, $\ad e:\g (i-2)\to\g (i)$ is injective
for $i\le 1$ and surjective for $i\ge 1$; 
%\par
%{\sf (\z_\g(e)} $(\ad e)^i : \g(-i)\rar\g(i)$ is one-to-one.
%{\rm (\z_\g(e)} $[\el,e]=\g(2)$.
%{\rm (v)} $\dim Z_G(e)=\dim \g(0)+\dim \g(1),\quad\dim(Z_G(e))^u=\dim
%\g(1)+\dim \g(2)$. 
\end{s}%
An element $e\in\N$ is said to be {\sl even\/}, if the eigenvalues of
$\ad h$ are even; 
$e\in\N$ is called {\sl distinguished\/}, if $\z_\g(e)(0)=\{0\}$. By a well-known
result of Bala-Carter and Vinberg, ``distinguished" implies 
``even".
The integer $\max\{i\in {\Bbb N}\mid \g(i)\ne 0\}$ is called the 
{\it height\/} of $e$ or of the orbit $G{\cdot}e$ and is denoted by $\hot(e)$. 
It is easily seen that $\hot(e)=m$ if and only if $(\ad e)^m\ne 0$ and
$(\ad e)^{m+1}= 0$. See 
\cite[Sect.\,2]{beyond} for some results concerning the height. 
\\[.6ex]
The notation related to the $\Bbb Z$-grading associated with a nilpotent 
orbit will be used throughout the paper. Whenever we discuss in the sequel
a $\Bbb Z$-grading associated with $e\in\N$, this means that $e$ is 
regarded as member of an $\tri$-triple and the grading is determined by $h$.
\\[.6ex]
{\bf Remark.} In the context of the adjoint representation (of $\g$)
we write $\z_\g(e)$ in place of $\g_e$. This algebra is also referred to
as {\it centraliser\/} of $e$, and not stabiliser. 
\\[.6ex]
Set $\de_\g(e)=\z_\g(\z_\g(e))$. Sometimes, this space is called
the {\it double centraliser} of $e$. It is also
easily seen that $\de_\g(e)$ is the centre of $\z_\g(e)$. 
The subalgebra $\de_\g(e)$ inherits the (nonnegative)
graded structure from $\z_\g(e)$.
R.\,Brylinski and B.\,Kostant have obtained several results on the 
graded structure of $\de_\g(e)$. For instance, they proved
that 
\begin{enumerate}
\item[$(\ast)$] \quad$\de_\g(e)(i)=\{0\}$ for $i\le 1$;
\item[$(\diamond)$] \quad $\dim\de_\g(e)(2)$ is equal to the number of
those simple ideals of $\g$ to which $e$ has nonzero projection.
\end{enumerate}
But whereas relation $(\ast)$ is proved rather
elementary (see~\cite[Prop.\,17]{rk1}), the proof of
$(\diamond)$ is very involved. 
In fact, the proof of Theorem~7.1 in~\cite{rk2} consists of a 
list of reference to previous results in two papers, which requires,
in particular, understanding some properties of a Poisson structure on
the algebra of regular functions on $G{\cdot}e$. 
%\,\cite{rk}
%and uses some results from Symplectic Geometry.
\\[.6ex]
We are going to give a simple self-contained proof of $(\diamond)$.
First, recall a required result of \cite{rk1}. Obviously,
the spaces $\de_\g(e)(2)$, $\de_\g(f)(-2)$, and $[f,\de_\g(e)(2)]$ have the
same dimension, say $d$. 
Set $\rr:=\de_\g(f)(-2)\oplus [f,\de_\g(e)(2)] \oplus\de_\g(e)(2)$.
%%Then $\ah\subset \p$ and $\dim\p=3d$.

\begin{s}{Proposition {\ququ \cite[Prop.\,20]{rk1}}} \label{rk} \\
The vector space $\rr$ is a subalgebra of $\g$. It is 
isomorphic to $(\tri)^d$.
\end{s}%
Set $\ah=\Bbbk e+\Bbbk h+\Bbbk f$. It is a subalgebra of $\rr$.
Let $\ah_1,\dots,\ah_d$ be the simple ideals in
$\rr$. Choose a basis $\{e_i,h_i,f_i\}$ 
of $\ah_i$ such that it is an $\tri$-triple, $e_i\in\de_\g(e)(2) $, and
$f_i\in \de_\g(f)(-2)$. It immediately follows from the definition of $\rr$
and \re{rk} that $\ah$ projects isomorphically to each $\ah_i$. We may 
therefore assume that $e=\sum_{i=1}^d e_i$. 
Let $\vp_i$ (resp. $\vp$) denote the fundamental weight of $\ah_i$
(resp. $\ah$) relative  to the Borel subalgebra 
$\Bbbk e_i+\Bbbk h_i$ (resp. $\Bbbk e+\Bbbk h$).
We shall use dominant weights in multiplicative notation to denote
irreducible $\rr$- and $\ah$-modules. For instance,
$\vp_1^{k_1}\dots\vp_d^{k_d}$ stands for the irreducible $\rr$-module
with highest weight $k_1\vp_1+\dots+k_d\vp_d$.

\begin{s}{Theorem {\ququ (R.\,Brylinski \& Kostant)}}    \label{d=1}  \\
Let $e$ be a nilpotent element in a semisimple Lie algebra $\g$. Then
$d=\dim\de_\g(e)(2)$ equals the number of those simple ideals of 
$\g$ to which $e$ has nonzero projection.
\end{s}%
\begin{proof} 1. Let $\g_1,\dots,\g_m$ be the simple ideals in $\g$
and $e=e^{(1)}+\dots+e^{(s)}$, where $0\ne e^{(j)}\in\g_j$ and $s\le m$. 
Since $\de_\g(e)(2)=\de_{\g_1}(e^{(1)})(2)\oplus\dots\oplus
\de_{\g_s}(e^{(s)})(2)$, it suffices to prove that $d=1$ whenever
$\g$ is simple. \par
2. Consider $\g$ as $\rr$-module: \\[.8ex]
\hbox to \textwidth{\refstepcounter{subsubsection}{\bf (\thesubsubsection)}
\label{first}
\hfil  
%\begin{equation} \label{first}
$\displaystyle \g\mid_\rr=\bigoplus_{k_1,\dots,k_d}m_{k_1\dots k_d}
\vp_1^{k_1}\dots\vp_d^{k_d}$.
\hfil} 
\vskip.6ex\noindent
%\end{equation}
Here $m_{k_1\dots k_d}$ is the multiplicity in $\g$ of the respective 
irreducible $\rr$-module.
Restricting further to $\ah$, we obtain:
\[  \g\mid_\ah=\bigoplus_{k_1,\dots,k_d}m_{k_1\dots k_d}
\vp^{k_1}\otimes\dots\otimes\vp^{k_d} 
\ .
\]
A highest weight vector in $\g$ relative to $\ah$ is nothing but a
homogeneous element of $\z_\g(e)$. Since each $e_i$ lies in $\de_\g(e)$, any 
highest weight vector in $\g$ relative to $\ah$ is also a highest weight
vector relative to $\ah_i$ and hence relative to $\rr$. This means that
the above
two decompositions have  ``the same" irreducible constituents. In other words,
$\vp^{k_1}\otimes\dots\otimes\vp^{k_d}$ is an irreducible $\ah$-module
for any $d$-tuple $(k_1,\dots,k_d)$ occurring in the 
\re{first}.
This clearly implies that any such $d$-tuple contains 
{\it at most one\/} nonzero coordinate. \par
%If $k_j$ is the unique nonzero coordinate of a $d$-tuple, we shall say that
%the respective summands in \re{first} are {\it nontrivial $\ah_j$-submodules}.
Let ${\frak L}_j$ ($j=1,\dots,d$) denote the sum of all nontrivial 
$\ah_j$-submodules, i.e., the sum of all summands in \re{first}, where
$k_j\ne 0$. Since each $d$-tuple contains
{\it at most\/} one nonzero coordinate,
 $\g=\z_\g(\rr)\oplus {\frak L}_1\oplus\dots\oplus {\frak L}_d$,
the sum of vector spaces. Note that $[{\frak L}_i,{\frak L}_j]=0$ ($i\ne j$), 
because any nonzero element would generate a $\rr$-submodule with a nontrivial 
action of {\it both\/} $\ah_i$ and $\ah_j$.
Denoting $\rr^{(j)}=\oplus_{i\ne j}\ah_i$, we see that 
${\frak L}_j\oplus \z_\g(\rr)=\z_\g(\rr^{(j)})$ is a reductive Lie subalgebra.
\par
%Now we are ready to make the last step. 
If $\h$ is a reductive subalgebra
of a simple Lie algebra $\g$, then the subalgebra $\es$ generated by
$\h^\perp$ equals $\g$.
Indeed, $\es$ is stable relative to both
$\ad(\h^\perp)$ and $\ad\h$, so that it is an ideal in $\g$,
see \cite[4.1]{kac}. 
Apply this to $\h=\z_\g(\rr^{(j)})$,
for arbitrary $j$. Here $\h^\perp=\bigoplus_{i\ne j}{\frak L}_i$, i.e.,
it is nonzero, if $d>1$. Then the relations of the previous paragraph
show that the
subalgebra generated by $\h^\perp$ lies in $\h^\perp\oplus\z_\g(\rr)$, i.e.,
it does not contain ${\frak L}_j$. Thus, $\g$ cannot be simple, if $d>1$.
\end{proof}%
%
%It was observed by P.\,Tauvel that this grading is always even, i.e.,
%$\de_\g(e)_i=0$ unless $i$ is even.
It follows from $\tri$-theory that
$\z_\g(e)\oplus [\g,f]=\g$, see Proposition~\ref{stab}(iii). 
It is equivalent to the fact that the 
restriction of the Killing form to $\z_\g(e)+\z_\g(f)$ is nondegenerate.
It is then rather natural to ask about a good complementary subspace to 
$\de_\g(e)$.
This question together with a conjectural answer was communicated to 
me by P.\,Tauvel in 1996. Here we give a proof for Tauvel's conjecture.

\begin{s}{Theorem}   \label{conj-P}
For any $\tri$-triple we have $\de_\g(e)\oplus [\g, \z_\g(f)]=\g$.
Equivalently, the Killing form is nondegenerate on
$\de_\g(e)\oplus\de_\g(f)\subset\g$.
\end{s}%
Before giving a proof, we state an auxiliary result.\\
Assume for a while that $\bbk=\Bbb C$.
Let $\theta$ be an involutory automorphism of $\g$.
It is a standard fact (see e.g. \cite[4.3.1]{t41})
that there exists an antilinear
involution $\sigma$ such that $\theta\sigma=\sigma\theta$ and
$\g^{ \theta\sigma}$ is a compact real form of $\g$.
We shall say that $\sigma$ is a {\it (complex) conjugation\/} associated with
$\theta$. In other words, $\theta$ is a Cartan involution for the real form
$\g^\sigma$.
%Then $\g^\sigma$ is a real form of $\g$ and
%$\g^\sigma\cap\g^\theta$ is a maximal compact subalgebra of $\g^\sigma$.
%An involutory  automorphism $\theta$ is said to be of {\it
%maximal rank\/} whenever
%its $(-1)$-eigenspace contains a Cartan subalgebra of $\g$.

\begin{s}{Proposition}  \label{complex}
Let $\theta\in\mathrm{Aut\,}\g$ be an involutory automorphism of $\g$
and $\sigma$ an associated complex conjugation. Suppose an $\tri$-triple
$\{e',h',f'\}$ satisfies the condition that $G{\cdot}h'$ has a nonempty
intersection with the $(-1)$-eigenspace of $\theta$. Then
$\{e',h',f'\}$ is $G$-conjugate to a triple $\{e,h,f\}$ such that
$\{e,h,f\}\subset\g^\sigma $, $\theta(h)=-h$, and $\theta(e)=-f$.
\end{s}\begin{proof}
This is  proved, although not stated explicitly, in \cite[Theorem\,1]{leva}.
\begin{comment} %%%%%%%%%%%%%%%%%%%%%%%
1. Let $\g_1$ denote the $(-1)$-eigenspace of $\theta$. Since
$\theta$ is of maximal rank, $G{\cdot}h\cap\g_1\ne\varnothing$.
... \\
2. So, we obtain $\{e,h,f\}\subset\g^\sigma$ and
$\theta({\Bbb R}e+{\Bbb R}h+{\Bbb R}f)\subset ({\Bbb R}e+{\Bbb R}h+{\Bbb R}f)$.
Note that the restriction of $\theta$ to 
${\Bbb R}e+{\Bbb R}h+{\Bbb R}f$ is not trivial, since 
$\g^{ \theta\sigma}$ is compact, whereas $\tri({\Bbb R})$ is not.
Thus, $\theta$ induces a nontrivial involution of $\tri({\Bbb R})$. It is 
not hard to verify that $\tri({\Bbb R})$ has a basis with the required 
properties with respect to $\theta$.
\end{comment} 
\end{proof}%

\noindent
{\sl Proof of \re{conj-P}.} 
1. First, we prove that the sum $\de_\g(e)+\de_\g(f)$ is direct. Indeed,
if $z\in \de_\g(e)\cap\de_\g(f)$, then it commutes with $f$ and $\z_\g(e)$.
Since the subalgebra generated by $f$ and $\z_\g(e)$ is the whole $\g$,
we see that $z=0$.
\par
2. Let $\theta$ be an involutory automorphism of {\it maximal rank\/},
i.e., its $(-1)$-eigenspace contains a Cartan subalgebra of $\g$.
Then the hypothesis of Proposition~\ref{complex} is satisfied for all
$\tri$-triples. We may therefore assume that $\{e,h,f\}$ are chosen
as prescribed in ~\re{complex}.
Then $\theta(\z_\g(e))=\z_\g(f)$ and hence
$\theta(\de_\g(e))=\de_\g(f)$, while $\sigma$ leaves all these spaces
intact. It follows that 
$\de_\g(e)\oplus\de_\g(f)$ is a $\theta$- and $\sigma$-stable
subspace of $\g$. Therefore it is the complexification of a subspace
in the compact real form $\g^{\theta\sigma}$. This shows 
that $\Phi$ is nondegenerate on 
$\de_\g(e)\oplus\de_\g(f)$. More explicitly, given $x\in\dge$,
we have $\sigma\theta(x)\in\de_\g(f)$ and $\Phi(x,\sigma\theta(x))\ne 0$.
\par
3. Using the invariance of the Killing form, one readily verifies that 
$[\g,\z_\g(f)]^\perp=\de_\g(f)$. This yields the equivalence stated in
the Theorem.
\qus%
We shall need a stronger result than $(\ast)$. It was proved by Tauvel in
\cite{patrice2}\footnote{Practically, this book seems to be
a virtual object because of the bankruptcy of the publisher.}. 
Here we give a shorter proof.

\begin{s}{Proposition}  \label{even}
For any $e\in\N$, we have  $\dge=\oplus_{i\ge 1}\dge(2i)$.
\end{s}\begin{proof}
Let $\h$ be a Cartan subalgebra of $\zge(0)$. Consider the reductive Levi
subalgebra $\el:=\z_\g(\h)$. It is clear that $\h$ is the centre of $\el$,
i.e., $\el=[\el,\el]\oplus\h$. By the construction of $\el$, we have
$\dge\subset\z_\el(e)$. Furthermore, $\dge\subset \z_{[\el,\el]}(e)$, since
$\dge(0)=0$. %%It remains to observe that 
Finally, the very definition of $\el$
implies that $e$ is distinguished in $[\el,\el]$ and therefore
$\z_{[\el,\el]}(e)$ is concentrated in positive even degrees.
\end{proof}%
{\bf Remark.}
The previous proof shows that the double centraliser of any $e\in\N$ is
isomorphic, as graded space, to the double centraliser of a distinguished
element in a smaller semisimple subalgebra. 
For this reason, all questions about  $\dge$ can be reduced to the
case, where $e$ is distinguished.
%\end{rem}

\sekt{On a conjecture of A.G.\,Elashvili
%On the index of centralisers in reductive Lie algebras
\label{ind=rk}}%
By Corollary~\ref{neravenstvo}, we have $\ind\q_\xi\ge\ind\q$ for any 
Lie algebra
$\q$ and any $\xi\in\q^*$. It is however easily seen that this inequality
can be strict.
\begin{rem}{Example} Let $\q$ be a Borel subalgebra of ${\frak sp}_4$.
It is not hard to verify that $\ind\q=0$ and there exists $\xi\in\q^*$
such that $\q_\xi$ is two-dimensional and Abelian, i.e., $\ind\q_\xi=2$.
\end{rem}%
The following remarkable conjecture claimes that such a phenomenon
cannot occur in reductive Lie algebras. For a reductive $\g$, we identify
$\g$ and $\g^*$ and consider centralisers $\z_\g(x)$ in place of stabilisers
$\g_\xi$ ($\xi\in\g^*$).

\begin{s}{Conjecture {\ququ (A.G.\,Elashvili)}}  \label{AG}
If $\g$ is a reductive Lie algebra, then $\ind\z_\g(x)=\rk\g$
for any $x\in\g$.
\end{s}%
Let us record several helpful (and easy) 
observations related to this conjecture.
\begin{itemize}
\item It suffices to prove the conjecture for the simple Lie algebras;
\item The conjecture is true, if $x$ is semisimple.
\item The Jordan decomposition in $\g$ shows that it is sufficient to
prove the conjecture only for the nilpotent elements. Indeed,
let $x=x_s+x_n$ be the Jordan decomposition, where $x_s$ is 
semisimple and $x_n$ is nilpotent.
Here $\el=\z_\g(x_s)$ is reductive, $\rk\g=\rk\el$,
$x_n$ is a nilpotent element of $\el$,
and the equality $\z_\g(x)=\z_\g(x_s)\cap\z_\g(x_n)$ means that
$\z_\g(x)=\z_\el(x_n)$.
\end{itemize}
For this reason, we shall only consider in the sequel the centralisers 
of {\sl nilpotent\/} elements.
It is worth noting that Problem 15 stated by R.Steinberg in his 
address to the I.C.M. in Moscow \cite[p.\,280]{icm66} (see also 
\cite[III.1.18]{ss}) can be regarded as a restricted 
version of this conjecture.
A solution to this problem is given by 
Kurtzke \cite{kurt} and Tauvel \cite{patrice}, independently. However,
both proofs exploit the classification of nilpotent elements and/or 
computer computations.
Here we give a short conceptual proof.

\begin{s}{Theorem}   \label{concept}
Let $e\in\g$ be a nilpotent element. Suppose $\z_\g(e)$ is
Abelian. Then $e$ is regular, i.e., $\dim\z_\g(e)=\rk\g$.
\end{s}\begin{proof}
Without loss of generality, we assume that $\g$ is simple.
Let $\{e,h,f\}$ be an $\tri$-triple. Recall from Section~\ref{RK-prosto}
the non-negative grading on $\z_\g(e)$ determined by $h$:
\[
         \z_\g(e)=\bigoplus_{i\ge 0}\z_\g(e)(i) \ .
\]
%As is well-known, $\z_\g(e)(0)$ is the centraliser of the $\tri$-triple.
Take any $z\in\z_\g(e)(0)$. By Proposition~\ref{stab}(ii), 
$z$ commutes with $f$ and, by the assumption,
with $\z_\g(e)$. From the $\tri$-theory it follows that the Lie subalgebra
generated by $f$ and $\z_\g(e)$ is $\g$. Hence $z$ lies in the centre of
$\g$, i.e., $z=0$. This means that $e$ is distinguished and therefore 
is even.
By Proposition~\ref{stab}(iii), we have 
\[
\dim\g(2i)-\dim\g(2i+2)=\dim\z_\g(e)(2i)=\dim\de_\g(e)(2i)
\]
for all 
$i\ge 0$. Then invoking Theorem~\ref{d=1}, we conclude that \\[.6ex]
%\begin{equation}  
\hbox to \textwidth{\refstepcounter{subsubsection}{\bf (\thesubsubsection)}
\label{codim1}
\hfil       $\dim\g(2)-\dim\g(4)=1$.         \hfil
}
\vskip.6ex\noindent
%%\end{equation}
This condition for a distinguished nilpotent element $x$ was considered
by T.A.\,Springer in \cite{tony}. 
Let $m=\hot(e)$
%%=\max\{i\in{\Bbb N}\mid \g(i)\ne 0\}$ 
and let $q$ be a nonzero element in
$\g(-m)$. Springer proved that $c:=q+e$ is a regular semisimple
element in $\g$ whenever \re{codim1} is satisfied, see Lemma~9.6 in 
[loc.\,cit.]
%We are going to prove that $\dim\z_\g(c)\ge\dim\z_\g(e)$, which completes the
%proof of the Theorem.
%\\[.7ex]
\begin{rems}{Claim}  \label{claim}
{\sl For every $x\in\z_\g(e)$ there exists a unique 
$\displaystyle z\in\bigoplus_{i<0}\g(i)$ such that $x+z\in\z_\g(c)$.}
\end{rems}%
As $e$ is distinguished, $x$ is nilpotent and therefore $x\notin\z_\g(c)$.
Therefore $z\ne 0$, if $x\ne 0$. The uniqueness of $z$ 
follows from the fact
that $\z_\g(c)$ contains no nonzero nilpotent elements. So, it remains to
prove that such a $z$ exists.
Consider the equation on $z$
\[
       0=[x+z,e+q]=[z,e]+[x,q] .
\]
(Note that $[q,z]=0$ in view of our assumption on $z$ and $q$.)
Since $\z_\g(e)$ is Abelian, $[x,q]$ is orthogonal
to $\z_\g(e)$ with respect to the Killing form on $\g$. Therefore 
$[x,q]\in \z_\g(e)^\perp=[\g,e]$, i.e. the above equation has a solution.
Because $[x,q]\in \oplus_{i\le 0}\g(i)$ and $e\in\g(2)$, we also deduce that
there exists a solution $z$ lying in $\oplus_{i<0}\g(i)$, as required.
\\[.6ex]
It follows from the Claim that 
$\dim\z_\g(e)\le\dim\z_\g(c)= \rk\g$, which completes the
proof of Theorem~\ref{concept}.
\end{proof}%
{\bf Remarks.} 1. In \cite{tony}, Springer gave a classification of all 
nilpotent orbits in simple Lie algebras satisfying \re{codim1}.
It turns that all such orbits in the classical Lie algebras
are necessarily regular, and there are a few nonregular nilpotent
orbits in the exceptional Lie algebras satisfying \re{codim1}, see Table~11
on p.196 in [loc.\,cit.]. 
However, the subregular orbit in type ${\Bbb G}_{2}$ 
is erroneously placed
in that list. Indeed, in this case $\dim\g(2)=4$ and
$\dim\g(4)=1$.

2. As was explained in \cite[III.4]{ss}, all results related to $\tri$-triples
that can be proved in characteristic zero
remain valid if $\textrm{char}(\bbk)\ge 4m+3$,
where $m$ is the height of the highest root. This means that
our proof of Theorem~\ref{concept} goes through under this constraint as well.
%%on the characteristic.
%
\begin{s}{Corollary}   \label{subreg}
Suppose $e\in\N$ is subregular (i.e., $\dim\z_\g(e)=\rk\g+2$).
Then $\ind\z_\g(e)=\rk\g$.
\end{s}\begin{proof}
From Corollary~\ref{neravenstvo} and Theorem~\ref{concept}, we conclude that
$\rk\g\le \ind\z_\g(e) < \rk\g+2$. Since $\dim\z_\g(e)-\ind\z_\g(e)$ is even,
we are done.
\end{proof}%
Our next goal is to show that Elashvili's Conjecture is true 
for `small' nilpotent orbits. Note that $\hot(e)\ge 2$ for all 
$x\in\N\setminus\{0\}$.
\begin{s}{Theorem}   \label{ht=2}
Suppose $\hot(e)=2$. Then $\ind\z_\g(e)=\rk\g$.
\end{s}\begin{proof}
Consider the ${\Bbb Z}$-grading of $\g$ associated with an $\tri$-triple 
$\{e,h,f\}$. We have $\g=\bigoplus_{i=-2}^{2}\g(i)$ and
$\z_\g(e)=\ka\oplus\g(1)\oplus\g(2)$, where $\ka$ is a reductive subalgebra
of $\g(0)$. Consider $f$ as element of $\z_\g(e)^*$, using the Killing form,
and compute its stabiliser $\z_\g(e)_f$. It readily follows from 
Proposition~\ref{stab} that $\z_\g(e)_f=\ka\oplus\g(2)$.
Combining Corollary~\ref{neravenstvo} and Theorem~\ref{thm:rais}
yields
\[
 \rk\g\le\ind\z_\g(e)\le\ind(\ka\oplus\g(2))=\ind(\ka:\g(2))+\ind \ka_\xi \ ,
\]
where $\xi$ is a regular element in the $\ka$-module $\g(2)^*$.
Thus, it suffices to prove that the last expression is equal to $\rk\g$. 
Note that $\g(2)$ is the isotropy representation for the 
{\sl affine\/} homogeneous 
space $G(0)/K$. This implies that $\g(2)\simeq\g(2)^*$ as $\ka$-modules
and there exists an open
subset $\Omega\in\g(2)$ such that the stabilisers $\ka_\xi$ ($\xi\in\Omega$)
are reductive and conjugate in $\ka$, see~\cite{luna}. 
Therefore we may assume that
$\ka_\xi$ is such a reductive generic stabiliser.
%%a generic stabiliser for the 
By \cite[3.3]{pMM}, $\ka$ is a `symmetric' subalgebra of $\g(0)$
whenever $\hot(e)\le 3$. In particular, $G(0)/K$ is a spherical homogeneous
space. In this case, one knows that
\[
  \dim \g(2)\md K =\rk G(0)-\rk \ka_\xi \quad (\xi\in\Omega) 
\]
%where $K_\xi$ is a generic stabiliser 
(see \cite[ch.II]{disser}). Here $\g(2)\md K=\spe (\bbk[\g(2)]^K)$.
Since $\g(2)$ is an orthogonal $K$-module, the field of invariants
$\bbk(\g(2))^K$ is the quotient field of $\bbk[\g(2)]^K$. Therefore
$\ind(\ka:\g(2))=\dim \g(2)\md K$ and, since $\ind \ka_\xi =\rk \ka_\xi$,
we conclude that $\ind(\ka:\g(2))+\ind \ka_\xi=\rk G(0)=\rk\g$.
\end{proof}%
Using some structure theory for the nilpotent orbits of height 3,
one can also show that Conjecture~\ref{AG} holds for such orbits.
This will be published elsewhere.
\begin{subs}{Current status of Elashvili's conjecture}
\end{subs}
The preceding exposition shows that the Conjecture is true for 
`large' and `small' nilpotent orbits. The area in between seems to be a
``terra incognita". However, explicit computations, performed by Elashvili
for all exceptional Lie algebras but
${\Bbb E}_8$, confirm the conjecture. He also succeeded
in proving the conjecture for ${\frak sl}_n$ (personal
communication). I hope that the details of his work will be published in 
the near future.

\sekt{The normaliser of the centraliser of a nilpotent element and its index
\label{normaliser}}%
In this section, $\g$ is a simple Lie algebra, and
we study the index of the Lie algebra 
$\n_\g(e):=\{s\in\g \mid [s, \z_\g(e)]\subset \z_\g(e) \}$, $e\in\N$.
The main structure results about $\n_\g(e)$ are obtained Brylinski and 
Kostant in \cite{rk1}.
It has to be mentioned that a nice exposition of these results, as well 
as some complements,  
is given in Tauvel's book, see \cite[ch.17]{patrice2}.
The following assertion contains, in a condensed form, all what we need
from that theory. 
\begin{s}{Theorem} %{\ququ (Brylinski \& Kostant)}} 
\label{ngx} 
Let $\{e,h,f\}$ be an $\tri$-triple.  Then \\ \indent 
{\sf (i)} \cite[Theorem 18]{rk1} \ 
$\n_\g(e)=\z_\g(e)\oplus [f,\de_\g(e)]$; in particular,
$\dim\n_\g(e)=\dim\z_\g(e)+\dim\de_\g(e)$. \\ \indent
{\sf (ii)} \cite[Theorem 23]{rk1} \ 
the Lie algebra $\q:=\n_\g(e)/\zge$ is solvable;
the image of the 1-dimensional subspace $\bbk h\subset [f,\de_\g(e)]$ 
is a maximal reductive subalgebra in it.  \\ \indent
{\sf (iii)} \cite[ch.17]{patrice2} \  
$\n_\g(e)=\{s\in\g \mid [s,e]\in \de_\g(e)\}$. 
\end{s}%
Let $N_\g(e)$ denote the connected algebraic
subgroup of $G$ with Lie algebra
$\n_\g(e)$. An immediate consequence of Theorem~\ref{ngx}(iii) 
is the claim that
\\[0.5ex]
\hspace*{.5pt}\refstepcounter{equation} {\bf (\theequation)}
\parbox{410pt}{\qquad $N_\g(e)$ has an open orbit in $\de_\g(e)$; namely,
$N_\g(e){\cdot}e$ is open in $\de_\g(e)$.
}\\[1ex]
We will be interested in the property that 
$N_\g(e)$ has an open orbit in $\de_\g(e)^*$. Recall that if an algebraic
group $A\subset GL(V)$ has an open orbit in $V$, then this does not imply
in general that $A$ has an open orbit in $V^*$. 
However, the following is true, see \cite{P} :
%\\[1ex]
%\hspace*{.5pt}\refstepcounter{equation} {\bf (\theequation)}
%\parbox{410pt}{\quad 
\centerline{\sl $A$ has finitely many orbits in $V$
if and only if it has finitely many orbits in $V^*$.
}
\\[.6ex]
From this we derive the following sufficient condition:
\\[1ex]
\hspace*{.5pt}\refstepcounter{equation} {\bf (\theequation)}
\parbox{410pt}{\quad If $N_\g(e)$ has finitely many orbits in $\de_\g(e)$,
then it has an open orbit in $\de_\g(e)^*$.
}
\begin{s}{Theorem}   \label{index-n}  \\
{\sf (i)} For any $e\in\N$, we have \par
 a) \ $\ind\n_\g(e)\ge \ind\z_\g(e)-\dim\de_\g(e)\ge \rk\g-\dim\dge$; \par
 b) \ $\ind(\n_\g(e),\z_\g(e)\ge \ind\zge-\dim\dge$;
\\
{\sf (ii)} 
suppose $N_\g(e)$ has an open orbit in $\de_\g(e)^*$. Then \\[.6ex]
\parbox{424pt}{\hfil
$\ind(\n_\g(e),\z_\g(e))=\ind\n_\g(e)=\ind\z_\g(e)-\dim\de_\g(e)$. \hfil} 
\end{s}\begin{proof}
Choose a basis $(e_1,\dots,e_n)$ for  $\z_\g(e)$ so that $(e_1,\dots,e_m)$
is a basis for $\de_\g(e)$, $m\le n$. Then 
$(e_1,\dots,e_n, [f,e_1],\dots, [f,e_m])$ is a basis for $\n_\g(e)$.
Write the matrix $\ck(\n_\g(e))$ in this basis. It has the following 
block structure: 
$\left(\begin{array}{ccc} 0   &  0    & {\frak D}     \\
                          0   &  {\frak C}    & \ast  \\  
                         -{\frak D}^t &  \ast & \ast
       \end{array}\right)$,
where  ${\frak D}$ (resp. ${\frak C}$) is a square matrix of order $m$ 
(resp. $n-m$).
Notice that some fragments of $\ck(\nge)$ are $\ck$-matrices in its own right.
Indeed, $\ck(\z_\g(e))=\left(\begin{array}{cc}  0   &  0     \\
                                            0   &  {\frak C}    
       \end{array}\right)$
and   $\ck(\n_\g(e),\de_\g(e))=(0\ 0\ {\frak D})$. It follows that 
$\rank {\frak C}=\dim\z_\g(e)-\ind\z_\g(e)$ and 
$\rank {\frak D}=
\dim\de_\g(e)-\ind(\n_\g(e),\de_\g(e))$, in view of~\re{ind(q,V)}.   

(i)a \ The block structure of $\ck(\nge)$ shows that 
\[
 \dim\nge -\ind\nge=\rank\ck(\nge)\le \rank {\frak C}+ 2\,\dim\dge \ ,
\]
which yields the first inequality. The second inequality stems from 
Corollary~\ref{neravenstvo}.

(i)b \ Applying Theorem~\ref{atilde} to $\tilde\q=\n_\g(e)$ and
$\q=\z_\g(e)$, we obtain \\[.6ex]
\hbox to \textwidth{\refstepcounter{subsubsection}{\bf (\thesubsubsection)}
\label{+}
\hfil     
$\ind\z_\g(e)+\ind\n_\g(e) \le \dim\de_\g(e) +2\,\ind(\n_\g(e),\z_\g(e))$.
\hfil
}
\vskip.6ex\noindent
Combining this with (a) yields the desired inequality.

(ii) The hypothesis means that $\ind(\n_\g(e),\de_\g(e))=0$, i.e.,
${\frak D}$ is non-singular. It also implies that
$\rank\ck(\n_\g(e))=2\,\rank {\frak D}+\rank {\frak C}$.
%%(in general one has only inequality $``\ge''$). 
Hence
\[
  \dim\n_\g(e) -\ind\n_\g(e) = 2\,\dim\de_\g(e) +\dim\z_\g(e) -
  \ind\z_\g(e)  \ ,
\]
which together with Theorem~\ref{ngx}(1) gives the second equality. 
To obtain  the first equality, notice that
$\ck(\n_\g(e),\z_\g(e))=\left(\begin{array}{ccc} 0   &  0    & {\frak D}     \\
                          0   &  {\frak C}    & \ast  \\  
            \end{array}\right)$
and therefore $\rank\ck(\n_\g(e),\z_\g(e))=\dim\de_\g(e) +\rank {\frak C}$. Thus,
\[
\dim\z_\g(e)-\ind(\n_\g(e),\de_\g(e))=
\dim\de_\g(e) +\dim\z_\g(e)-\ind\z_\g(e) \ ,
\]
which completes the proof.
\end{proof}%
{\bf Remark.} 
The validity of Conjecture~\ref{AG} would imply that
$\ind\n_\g(e)=\rk\g-\dim\de_\g(e)$  under the hypothesis of 
Theorem~\ref{index-n}(ii).
\\[1ex]
It is difficult to directly verify whether $N_\g(e)$ has an open orbit in
$\dge^*$. However, it appears to be easy to describe a wide
class of examples, where $N_\g(e)$ has finitely many orbits in $\dge$.
\\
Because $h$ has even eigenvalues on $\dge$, it is convenient to replace it
by $\tilde h=\frac{1}{2}h$.  Let $1=m_1< m_2\le m_3\le\dots\le m_l$ be
the $\tilde h$-eigenvalues on $\dge$, counted with multiplicities, i.e.,
$l=\dim\dge$. (That $m_2>1$ follows from Theorem~\ref{d=1}.)
Let $(e_1,\dots,e_l)$ be a respective basis of $\tilde h$-eigenvectors.
Notice that $e_1=e$. 
\begin{s}{Proposition}   \label{progress}  The group
$N_\g(e)$ has finitely many orbits in $\dge$ if and only if \\[.6ex]
\hbox to \textwidth{\ $(\heartsuit_1)$ \hfil 
$[[f,e_j],e_i]=\ap e_{i+j-1}$ for all $i,j$ such that $i+j-1\le l$\hfil}
\vskip.7ex\noindent
and for some
$\ap=\ap_{ij}\in\bbk\setminus\{0\}$.
In particular, \\[.6ex]
\hbox to \textwidth{\ $(\heartsuit_2)$ \hfil 
$m_i+m_j-1=m_{i+j-1}$  \hfil}
\noindent
and
$\tilde h$ has a simple spectrum on $\dge$.
\end{s}\begin{proof*}
Set $Z_G(e)^o=Z_G(e)\cap N_\g(e)$. It is the
uneffectivity kernel for the $N_\g(e)$-action  on $\dge$. Letting
$Q=N_\g(e)/Z_G(e)^o$, we conclude from Theorem~\ref{ngx}(ii) that 
%%$Q$ is a (connected) solvable algebraic group and 
$Q\simeq T_1\,\semidir\, Q^{u}$,
where $T_1$ is a 1-dimensional torus and $Q^{u}$ is the unipotent radical
of $Q$.
Of course, we choose $T_1$ so that all $e_i$'s are eigenvectors of it.
Considering the $Q$-action on $\dge$, we may identify $\q$, the Lie algebra of
$Q$, with the space $[f,\dge]$.
Letting $L_i:=\oplus_{j\ge i}\bbk e_j$, we obtain a $Q$-stable complete 
flag $\dge=L_1\supset L_2\supset\dots\supset L_l$. 

a) Suppose conditions $(\heartsuit_1)$ are satisfied. Obviously, these 
conditions mean that $[\q, e_i]=L_i$ for all $i$. In other words,
$Q{\cdot}e_i$ is dense in $L_i$. More precisely, we have 
$Q^{u}{\cdot}e_i=e_i+L_{i+1}$ and $T_1{\cdot}e_i=(\bbk\setminus\{0\})e_i$.
Thus, $L_i=Q{\cdot}e_i\cup L_{i+1}$, and finiteness follows.

b) Conversely, suppose $Q$ has finitely many orbits in $\dge$.
Then each $L_i$ must contain a dense orbit. This already implies that
all $m_i$'s to be different. Indeed, assume that
%It follows from \ref{ngx}(i),(ii) that
%$Q{\cdot}e=\{ \sum_i x_ie_i\mid x_1\ne 0\}$. That is, the complement of 
%the open orbit is the hyperplane $\oplus_{i\ge 2}\bbk e_i$. Then we argue
%`inductively'. If 
$m_i=m_{i+1}$ for some $i$. Let $x=\sum_{j\ge i}x_je_j$ be a generic element of 
$L_i$. Then $x\mapsto x_i/x_{i+1}$ is a $Q$-invariant rational function on 
$L_i$, and therefore a dense orbit cannot exist. Thus $m_i < m_{i+1}$ for
all $i$. Next, let $x\in L_i$ be such that $[\q, x]=L_i$. Clearly
$[[f,e_j],x]=0$ whenever $j>l+1-i$. Writing the condition that the vectors
$[[f,e_j],x]$  ($j\le l+1-i$) are linearly independent, we obtain exactly 
conditions  $(\heartsuit_1)$ with fixed $i$.
%$Q$ cannot have an open orbit in
%$\dge$. Indeed, for any $q\in Q^{u}$ we have $q{\cdot}e_j-e_j\in L_{i+2}$
%($j=i,i+1$); whereas the toral part of $Q$ is 1-dimensional. Therefore
\end{proof*}%
\begin{s}{Corollary {\ququ (of the proof)}}   \label{chislo}  
If $N_\g(e)$ has finitely many orbits in $\dge$, then the number of orbits
is $\dim\dge+1$ and each orbit closure is an affine space. \qu
\end{s}%
Note that conditions $(\heartsuit_2)$ are vacuous for $i=1$, whereas for $i=2$,
we obtain that $\{m_i\}$ is an arithmetic progression with first term 1.
For $i>2$, condition $(\heartsuit_2)$ yields no new constraints on $\{m_i\}$.
\begin{s}{Theorem}   \label{3serii}  \\
Condition \ref{progress}$(\heartsuit_1)$ is satisfied in the
following cases:

{\sf (i)}\quad for every nilpotent element, if 
$\g\in \{ {\frak sl}_{n}, {\frak sp}_{2n}, {\frak so}_{2n+1}\}$;

{\sf (ii)}\quad for any nilpotent element whose partition contains at least
three parts, if $\g={\frak so}_{2n}$.

{\sf (iii)}\quad if\/ $\dim\dge\le 2$.
\\ 
In particular, if $e\in\N$ falls in either of these cases, then 
\par
  a)\quad $N_\g(e)$ has finitely many orbits in $\dge$; \par
  b)\quad $\ind\nge=\ind(\nge,\zge)=\ind\zge-\dim\dge$. 
\end{s}\begin{proof}
Let $e\in {\frak sl}_n$ be a nilpotent matrix and let $\{e,h,f\}$ be an 
$\tri$-triple. Then, denoting by $e^i$ the usual $i^{th}$ matrix power,
one easily proves by induction that
\[
[e^i,f]=e^if-fe^i=\displaystyle
\sum_{\ap+\beta=i-1} e^\ap h e^\beta \textrm{\quad and \quad }
[h, e^j]=2j{\cdot}e^j \ .
\]
Then
\[
   [[e^i,f],e^j]=\sum_{\ap+\beta=i-1}e^\ap h e^{\beta +j} -
\sum_{\ap+\beta=i-1}e^{\ap+j} h e^{\beta}=
\sum_{\ap+\beta=i-1}e^\ap [h,e^j] e^\beta=
   2ij{\cdot}e^{i+j-1} \ .
\]   
This equality is an incarnation of \ref{progress}$(\heartsuit_1)$ for
classical Lie algebras. 
Indeed, it is well known (see e.g. \cite{kurt1})
that, for $e\in {\frak sl}_{n}$, $\dge$ is the span of
all nonzero powers $e^j$. For ${\frak sp}_{2n}$ and 
${\frak so}_{2n+1}$, the same holds with the 
{\sl odd\/} powers of $e$. This proves (i).
For ${\frak so}_{2n}$, it is known that $\dge$ is the span of all odd powers
if and only if the partition of $e$
has at least three parts, whence (ii). Part (iii) is obvious.\\
The conclusions a) and b) follows by \ref{index-n}(ii) and \ref{progress}.
\end{proof}%
{\bf Example.} Take a subregular nilpotent element $e\in {\frak so}_8$.
Here $\dim\zge=6$, $\dim\dge=3$, and the corresponding partition is
$(5,3)$, i.e., it has two parts. It is easy to find that
here $m_1=1, m_2=m_3=3$. Thus, condition
\ref{progress}$(\heartsuit_2)$ is not satisfied for $i=j=2$. Therefore 
$N_\g(e)$ (or $Q$) has infinitely many orbits in $\dge$.
This shows that assumptions \ref{3serii}(ii),(iii) cannot be weakened
in general. The fact that $m_2=m_3$ also implies that $N_\g(e)$ does not have
an open orbit in $\dge^*$.

\sekt{The normaliser of the centraliser of a regular
nilpotent element
\label{principal}}%
As above, $\g$ is simple. In this section, $e\in\N$ is a 
\un{re}g\un{ular}
(= principal) nilpotent element. In this case, an $\tri$-triple
$\{e,h,f\}$ is said to be {\sl principal\/}. It is worth noting that
principal $\tri$-triples were first considered by E.B.\,Dynkin in
\cite{ebd50}. 
Our aim in this section is to prove that $\nge$ is a Frobenius Lie algebra.
Note that, for ${\frak sl}_n, {\frak so}_{2n+1}, {\frak sp}_{2n}$, and 
${\Bbb G}_2$, this already follows from
Theorem~\ref{3serii}. But the proof to be given below is absolutely general and
does not appeal to classification. Along the way we prove several properties of 
principal $\tri$-triples, which seem to be new. 
\\
We keep the notation of Section~\ref{RK-prosto}.
Recall that $p=\rk\g$ and $\{\ap_1,\dots,\ap_p\}$ are the simple roots.
Without loss of generality, we assume that $h=h_+$ (see Section~\ref{RK-prosto}).
Then $\g(0)=\te$ and 
%%$\oplus_{i\ge 1} \g(2i)=\ut_+$. Actually,
$\g(2i)=\oplus_{\hot(\ap)=i}\g_\ap$, where $\hot(\ap)$ is the usual height
of a root $\ap$. In particular, $\g(2)=\oplus_{\ap\in\Pi}\g_\ap$.
Recall some standard invariant~theoretic features of this situation, see
\cite{ko59},\,\cite{ko63}.
\begin{itemize}
\item The eigenvalues of $h$ on $\zge$ are $\{2m_1,\dots, 2m_p\}$,
where $\{m_1,\dots, m_p\}$ are the {\it exponents\/} of (the Weyl group of)
$\g$. We order the exponents so that $1=m_1 < m_2 \le \dots \le m_p$;
\item $\bbk[\g]^G$ is a graded polynomial algebra. The degrees of
free generators are $\{m_i+1\}_{i=1,\dots,p}$. The quotient map
$\pi_G : \g \to \g\md G$ is flat and $\pi_G^{-1}(\pi_G(0))=\N$;
\item Let $M$ be any $\ad h$-stable complement of $[\g,e]$ in $\g$. 
Then the eigenvalues of $h$ on $M$ are $\{-2m_1,\dots, -2m_p\}$.
The affine subspace $e+M$ is transversal to $G{\cdot}e$ 
at $e$. The restriction homomorphism
$\bbk[\g]\to \bbk[e+M]$ induces an isomorphism
$\bbk[\g]^G\isom \bbk[e+M]$. This
means, in particular, that $(e+M)\cap\N=\{e\}$. Such an affine subspace 
is said to be
%%We shall say the $e+M$ is 
a {\it special transversal slice\/} (at $e$).
\end{itemize}
The usual choice for $M$ is $\z_\g(f)$; however we shall see that there are
some other natural possibilities. 
\\
For each $\ap\in\Delta$, define $h_\ap\in\te$ by the formula $\Phi(h_\ap,x)=
\ap(x)$ for all $x\in\te$.
Let $e_\ap$ be a nonzero element in $\g_\ap$.
We normalise them so that $\Phi(e_\ap,e_{-\ap})=1$. Then $[e_\ap,e_{-\ap}]=
h_\ap$.
We have $\hot(e)=2m_p$ and
%%the ``highest'' nonzero space 
$\g(2m_p)=\g_\lb$, where
%%The corresponding root of $\g$ is $-\lb$, where 
$\lb$ is the highest root in
$\Delta_+$. 
%We will assume that $\Phi(e_\lb,e_{-\lb})=1$.
%$\{e_\lb, [e_\lb, e_{-\lb}], e_{-\lb} \}$ is an $\tri$-triple.
Since $\dim\g(2)-\dim\g(4)=1$, Springer's result \cite[9.6]{tony} says that
$c=e+e_{-\lb}$ is regular semisimple, cf.
the proof of Theorem~\ref{concept}.
Applying Claim~\ref{claim} in this situation, we see that, for $e_\lb\in\zge$,
%us to obtain a complete description
%of the centraliser of $c$. In particular, as $e_\lb\in\zge$,
there exists a unique $\hat f\in \ut_-$ such that $e_\lb+\hat f\in\z_\g(c)$.
\begin{s}{Proposition}    \label{hat-f}
{\sf (i)} $\hat f$ lies in $\g(-2)$; it is a regular nilpotent 
element.\\ \hspace*{3.8cm}
{\sf (ii)} $\z_\g(c)\subset \zge\oplus\zghf$.
\end{s}\begin{proof*}
(i) The element $\hat f\in\ut_-$ is being determined by the relation
\[
   [e,\hat f]=[e_\lb, e_{-\lb}]=h_\lb\in\g(0) \ .
\]
This already shows that $\hat f\in\g(-2)$.
Letting $e=\sum_{i=1}^p X_ie_{\ap_i}$ and 
$\hat f=\sum_{i=1}^p Y_ie_{-\ap_{i}}$, we obtain 
$[e,\hat f]=\sum_{i=1}^p (X_iY_i) h_{\ap_i}$. On the other hand, 
$h_\lb=\sum_{i=1}^p n_ih_{\ap_i}$, where $\{n_i\}$ are defined by the
equality $\lb=\sum_{i=1}^p n_i\ap_i$. Hence all $n_i \ne 0$,
Since $e$ is regular, we have all $X_i\ne 0$. It follows that all
$Y_i\ne 0$, too. Thus, $\hat f$ is regular as well.
%we see that
%$\hat f$ must have nonzero projection to all root spaces $\g_{-\ap_i}$.
%Therefore $\hat f$ is regular as well. 
\\
(ii) Claim~\ref{claim} applied in our situation gives the following: For any
$x\in\zge$ there exists a unique $z\in\ut_-$ such that $x+z\in\zgc$. Since
$\dim\zge=\dim\zgc$, this construction gives the whole centraliser of $c$.
Thus, each $y\in\zgc$ has only `positive' and `negative' parts;
$y=y_++y_-$, where $y_+\in\ut_+$ and $y_-\in\ut_-$. The totality of all positive
parts forms $\zge$. Because  $\zgc$ is Abelian, the totality of all negative 
parts forms an Abelian subalgebra of dimension $p$. Since $\hat f$ is one
of such negative parts and it is regular, the totality of all negative parts 
is precisely $\zghf$.
\end{proof*}%
\begin{rems}{Remark}
Because $\dim [\g_\lb, \g_{-\lb}]=1$, the line $\bbk\hat f$ and the centraliser
$\zghf$ do not depend on the particular choice of $e_\lb$ and $e_{-\lb}$.
It is easily seen that $\hat f$ is proportional to $f$
if and only if $h_\lb$ is proportional to $h$ if and only
if $\g={\frak sl}_2$ or ${\frak sl}_3$.
That is, in general, these two elements are completely different. Nevertheless,
since $\hat f$ is regular and $\hat f\in\g(-2)$, the space 
$\zghf$ is $\ad h$-stable and the eigenvalues of
$h$ on $\zghf$ are the same as on $\z_\g(f)$.
\end{rems}%
Until Theorem~\ref{slice2} we will work with $\bbk=\Bbb C$. 
Let $\theta$ be an involution of $\g$ of maximal rank that acts as \ $-id$ 
on $\te=\g(0)$. Let $\sigma$ be an associated complex conjugation. 
Then $\theta(\g_\ap)=\g_{-\ap}$  and $\sigma(\g_\ap)=\g_{\ap}$
for each $\ap\in\Delta$. Note that $\g^\sigma$ is a split real form of $\g$
and all $h_\ap$ ($\ap\in\Delta$) lie in $\g^\sigma$.
In view of Proposition~\ref{complex},
we may assume that $e,h,f\in \g^\sigma$ and $\theta(e)=-f$.
Set $\te_{\bbb R}=\te\cap\g^\sigma$. Then $\te=\te_{\bbb R}\oplus
\sqrt{-1}\,\te_{\bbb R}$ and $\sqrt{-1}\,\te_{\bbb R}\subset\g^{\sigma\theta}$.
Set $T_{\bbb R}:=exp(\te_{\bbb R})$ and $T_c=exp(\sqrt{-1}\,\te_{\bbb R})$.
Then $T=T_{\bbb R}{\cdot}T_c$ is a maximal torus in $G$.
\begin{s}{Lemma}  \label{nevyr}
Suppose $t\in T_{\bbb R}$.  %%:=exp(\te_{\bbb R})\subset G_{\bbb R}$. 
Then $\Phi$ is nondegenerate
on $\zge\oplus t{\cdot}\z_\g(f)$.
\end{s}\begin{proof}
Take any $s\in T$ such that $s^2=t$. Writing $s=s_{\bbb R}s_c$ (in the
obvious notation), we obtain $(s_c)^2=1\in T$. Set
$V=\zge\oplus t{\cdot}\z_\g(f)$. Then $s^{-1}{\cdot}V=
s^{-1}{\cdot}\zge\oplus s{\cdot}\z_\g(f)$.
Since 
%It follows from the previous constructions that
$\sigma\theta(s^{-1}{\cdot}\zge)=(s_c^{-2} s){\cdot}\z_\g(f)=
s{\cdot}\z_\g(f)$, the space $s^{-1}{\cdot}V$ is
the complexification of a space in the compact real form 
$\g^{\sigma\theta}$. Thus, $\Phi$ is nondegenerate on  $s^{-1}{\cdot}V$ 
and hence on $V$.
\end{proof}%
The following result says that $\hat f$ can be used in place of
$f$ for constructing a special transversal slice to $G{\cdot}e$ at $e$.
\begin{s}{Proposition}    \label{slice1}
$\z_\g(\hat f)\oplus [\g,e]=\g$, or equivalently, $\Phi$
is nondegenerate on $\zge\oplus\zghf$. 
\end{s}\begin{proof}
Since both $f$ and $\hat f$ are regular nilpotent elements in $\g(-2)$,
there exists a unique $t\in T$ such that $t{\cdot}f=\hat f$.
To order to define $\hat f$, we used the relation $[e,\hat f]=h_\lb$.
Because $\h_\lb\in\g^\sigma$, we have $\hat f\in \g^\sigma$ as well, and 
therefore $t\in T_{\bbb R}$. Since $\zghf=t{\cdot}\z_\g(f)$,
we conclude by Lemma~\ref{nevyr}.
%%has a square root, say $t$, which lies again in.
\end{proof}%
Yet another possibility for a special transversal slice is provided 
by the next theorem.
\begin{s}{Theorem}    \label{slice2}
$[e_{-\lb},[f,\zge]]\oplus [\g,e]=\g$.
\end{s}\begin{proof}
Set $\ca= [e_{-\lb},[f,\zge]]$ and $\cb=[e_{-\lb},\zge]$. 
We begin with proving that $\ca$ has the required dimension,
i.e., $p$. Since
$c=e_{-\lb}+e$ is regular semisimple, $\z_\g(e_{-\lb})\cap\zge=\{0\}$.
Hence $\dim\cb=p$. Because $\zge$ is Abelian, we see that $\cb$ is orthogonal
to $\zge$ with respect to the Killing form. In other words, $\cb\subset [\g,e]$
and hence $\cb\cap\z_\g(f)=\{0\}$. Consequently, $\dim\ca=p$, as
$\ca=[f,\cb]$.
\\
Using the fact the Cartan subalgebra $\zgc$ is ``diagonally'' embedded
in $\zge\oplus\zghf$, we obtain
\[
\cb=[e_{-\lb}+e,\zge]=[e_{-\lb}+e,\zghf]=[e,\zghf] \ .
\]
It follows that $\ca=[f,[e,\zghf]]$. Given $x\in\zghf$, we have
$[f,[e,x]]=-[h,x]+[e,[f,x]]$. Since $\z_\g(h)\cap\zghf=\{0\}$,
each nonzero element of $\ca$ is a sum of a nonzero element
in $\zghf$ and an element in $[\g,e]$. Thus, the validity of 
Theorem~\ref{slice2}
follows from that of Proposition~\ref{slice1}.
\end{proof}%
%
%{\bf Remark.} 
%It is worth noting that the statement of Theorem~\ref{slice2} does not 
%exploit the construction of $\hat f$ and various assumptions related to
%real forms of $\g$. The description of a special transversal slice is given 
%entirely in terms of the $\tri$-triple.
%\\
The following is our main result in this section.
\begin{s}{Theorem}    \label{frobenius}
$\nge$ is a Frobenius Lie algebra; more precisely,
$e_{-\lb}$, considered as element of $\nge^*$, is a representative of
the open $N_\g(e)$-orbit.
\end{s}\begin{proof}
Using the Killing form, each $x\in\g$ can be considered as element of
the dual space $\nge^*$. 
\\
Let us describe a convenient model for  $\nge^*$.  In our situation,
$\nge=\zge\oplus [f,\zge]=\Ker (\ad e)^2$. It follows that $\nge^\perp
=\Ima(\ad e)^2$. Since $\g=\Ima(\ad e)^2\oplus\Ker (\ad f)^2$ (which
actually holds for any distinguished element), we may identify $\nge^*$
with $\Ker (\ad f)^2$. The $\nge$-module structure on 
$\Ker (\ad f)^2$ is given by the usual Lie bracket in $\g$ coupled with
the subsequent projection to $\Ker (\ad f)^2$, with kernel $\Ima(\ad e)^2$.

Using the notation of the proof of Theorem~\ref{slice2}, we have
$[e_{-\lb},\nge]=\ca+\cb$. It is already proved therein
that $\dim\ca=\dim\cb=p$, $\cb\subset [\g,e]$, and 
%whereas Theorem~\ref{slice2} says that
$\ca\cap [\g,e]=\{0\}$. Hence $\dim(\ca+\cb)=2p$. Thus, we only need to show
that $(\ca+\cb)\cap \Ima(\ad e)^2=\{0\}$. As was noted in the proof of
Theorem~\ref{slice2}, a stronger condition holds for$\ca$; 
namely, $\ca\cap \Ima(\ad e)=\{0\}$. Next, we can restate 
Theorem~\ref{slice2} as follows: $\Phi$ induces a nondegenerate 
pairing between $\ca=[f,\cb]$ and $\zge=[\g,e]^\perp$. Using the invariance of
$\Phi$, this also yields a nondegenerate 
pairing between $\cb$ and $[f,\zge]$. It follows that 
$\cb\cap [f,\zge]^\perp=\{0\}$.
Since $[f,\zge]\subset\nge$, we see that $\cb\cap\nge^\perp=\{0\}$.
This means that the image of $[e_{-\lb},\nge]$ in $\g/\nge^\perp=
\nge^*$ is of dimension $2p$, which
completes the proof. % of the theorem.
\end{proof}% 
Making use of the matrix approach to the index, we can deduce from
Theorem~\ref{frobenius} that a certain matrix is nonsingular.

\begin{s}{Corollary}    \label{matrix}
Let $\{e_1,\dots,e_p\}$ be a basis for $\zge$. Then 

{\sf (i)} the $p\times p$ matrix
${\frak D}=\bigl([e_i,[e_j,f]]\bigr)_{i,j=1}^p$ is nonsingular and symmetric. 
(The entries of ${\frak D}$, which belong to $\zge$, are 
regarded as elements of the symmetric algebra ${\mathcal S}(\zge)$); 

{\sf (ii)}  up to a constant multiple, $\textrm{\rm det}({\frak D})
=(e_\lb)^p$. In particular, 
the matrix ${\frak D}(e_{-\lb})=
\bigl(\Phi([e_i,[e_j,f]],e_{-\lb})\bigr)_{i,j=1}^p$ is nonsingular.
\end{s}\begin{proof}
The symmetricity of ${\frak D}$ is obvious.
By Theorem~\ref{ngx}(i), $(e_1,\dots,e_p,[e_1,f],\dots,[e_p,f])$ is a basis
for $\nge$. The matrix $\ck(\nge)$ in this basis is
$\left(\begin{array}{cc}      0     &  {\frak D} \\
                       -{\frak D}^t &  \ast 

\end{array} \right)$. This proves (i). \\
%%We also know from 
By Theorem~\ref{frobenius}, %that
$\ck(\nge)$ remains nonsingular after having been evaluated against
$e_{-\lb}\in\nge^*$. Hence ${\frak D}(e_{-\lb})$ is nonsingular. Furthermore,
$\textrm{\rm det}({\frak D})$ vanishes precisely on the complement of the open
$N_\g(e)$-orbit in $\nge$. In our situation, $e_{-\lb}\in\nge^*$ is the 
unique weight vector of minimal weight with respect to $h$. Let 
$W\subset \nge^*$ be {\it the\/}
$\ad h$-stable complementary hyperplane. Since $N_\g(e)$ is a semi-direct
product of a 1-dimensional torus and a unipotent group,
it is easily seen that
$\nge^*\setminus N_\g(e){\cdot}e_{-\lb}=W$ (cf. proof of Theorem~\ref{index-n}).
Because $W$ is the annihilator of $e_\lb$, we are done.
\end{proof}% 
More specifically, let us choose a basis for $\zge$ so that
$e_i\in \g(2m_i)$. Then $d_{ij}:=[e_i,[e_j,f]]\in \zge(2m_i+2m_j-2)$.
Therefore it is 0 whenever $m_i+m_j-1\ > m_p$.
Recall that the exponents satisfy the property that
$m_i+m_{p+1-i}=m_p+1$ for all $i$ (this sum is called the 
{\it Coxeter number\/} of $\g$). Consequently, if the exponents are different,
then $d_{ij}=0$ for $i+j> p+1$. That is, all entries below the antidiagonal 
are equal to 0.
% while the entries on the antidiagonal are equal to
%$e_p=e_{\lb}$ up to a constant multiple. 
Since ${\frak D}$
is nonsingular, the antidiagonal entries must be non-zero. Thus, we derive
by `grading' reasons the following.

\begin{s}{Corollary}    \label{matrix2}
Suppose all the exponents of $\g$ are different, and let $e_i\in\zge(2m_i)$.
Then $[e_i,[e_{p+1-i},f]]=\mu_i e_{\lb}$ with some 
$\mu_i\in \bbk\setminus\{0\}$.
\end{s}%
For $\g={\frak so}_{4n}$, there are two equal exponents. However, it is also
possible to choose basis vectors in the corresponding 2-dimensional space so 
that the whole ${\frak D}$ to be triangular.

\sekt{Some open problems\label{voprosy}}%
In this section, $x\in\g$ is an arbitrary nilpotent element.
It is proved in Section~\ref{normaliser} that, under certain assumptions,
$\ind\n_\g(x)=\ind\z_\g(x)-\dim\de_\g(x)$. 
Since those assumptions are satisfied for
`almost all' classical cases, it is tempting to assume that this phenomenon
is completely general. Combining this with Elashvili's conjecture,
we arrive at the following nice statement.

\begin{s}{Conjecture}    \label{c1}
Let $x\in\N\subset\g$. Then $\ind\n_\g(x)=\rk \g-\dim\de_\g(x)$.
\end{s}%
If true, this conjecture would give an explanation to the 
fact observed by Kurtzke \cite{kurt} that $\dim\de_\g(x)\le \rk\g$ for all 
$x$. It is worth noting that 
$\dim\de_\g(x)=\rk\g$ if and only if either $x$ is regular of
$\g={\Bbb G}_2$ and $x$ is subregular. So, the
conjecture implies that $\n_\g(x)$ is Frobenius only in these two cases.
\\
Actually, I believe that the following stronger assertion is true.
\begin{s}{Conjecture}    \label{c2}
%%Let $x\in\N\subset\g$. Then 
$\ind(\n_\g(x),\z_\g(x))=\rk \g-\dim\de_\g(x)$.
\end{s}%
Indeed, it follows from Theorem~\ref{index-n}(i) and \re{+}
that Conjecture~\ref{c2} implies \re{c1}. 
\\[.6ex]
The previous exposition shows that the last conjecture is valid
in two extreme cases: \par
a) \ $x$ is regular (by Theorem~\ref{frobenius}); \par
b) \ $\hot(x)=2$  (here $\de_\g(x)=\bbk x$, by
Theorem~\ref{d=1}. Then combine Theorem~\ref{ht=2} and 
Theorem~\ref{index-n}(ii)).
\\
A proof of the last conjecture may consists of two steps. 
The first step is to prove Elashvili's conjecture, and the second is to prove
the equality
$\ind(\n_\g(x),\z_\g(x))=\ind\z_\g(x)-\dim\de_\g(x)$. 
The latter has a transparent 
geometric meaning. By dimension reason, we have the inequality
\[  
  \max_{\eta\in\z_\g(x)^*}\dim\z_\g(x){\cdot}\eta+\dim\de_\g(x)
  \ge \max_{\xi\in\z_{\frak g}(x)^*}\dim\n_\g(x){\cdot}\xi \ ,
\] 
and the equality in the second step is the same as equality here.
%To attack these conjectures in general, one may try to find out a suitable 
%inductive procedure.
\\[.7ex]
$\bullet$ \ 
It would be interesting to find proofs for Propositions~\ref{conj-P}
and \ref{slice1} that do not appeal to $\Bbb C$. Note  
in this regard that 
Lemma~\ref{nevyr} is wrong for arbitrary $t\in T$.
\\[.7ex]
$\bullet$ \ 
In conclusion, I state a problem related to parabolic subalgebras. 
Corollary~\ref{levi} says that $\ind\p+\ind\p^u\le\dim\el$. It is likely
that there is also a lower bound, namely,
$\ind\p+\ind\p^u\ge\rk\g$. 
%I hope to discuss this subject in a forthcoming 
%publication.
%\noindent
%${\frak abcdefghijklmnopqrstuwxyz}$ \\
%${\frak ABCDEFGHIJKLMNOPQRSTUWXYZ}$

%
%\vno{2.5} \indent
%{\footnotesize 
%\parbox{195pt}{% 
%{\it ul. Akad. Anokhina, d.30, kor.1, kv.7 \\ 
%Moscow 117602 \quad Russia} \\ {\tt panyush@mccme.ru }}}
%

\end{document}